%% file: symcov.tex
\documentclass{article}
\usepackage{amsmath,amssymb}
\usepackage{graphicx}

\newtheorem{theorem}{Theorem}
\newtheorem{definition}{Definition}
\newtheorem{lemma}[theorem]{Lemma}
\newtheorem{rem}[theorem]{Remark}

\newtheorem{alg}{Algorithm}

\newcommand{\inte }{\mathrm{int}\,}

\newcommand{\sgn }{{\rm sgn}\,}
\newcommand{\bd }{\partial}

\newcommand{\dom }{\mathrm{dom} \,}
\newcommand{\dir }{\mathrm{dir} \,}
\newcommand{\Fix }{\mathrm{Fix} \,}

\newcommand{\cover}[1]{\stackrel{#1}{\Longrightarrow}}

\newcommand{\bicover}[1]{\stackrel{#1}{\Longleftrightarrow}}
\newcommand{\invcover}[1]{\stackrel{#1}{\Longleftarrow}}

\def\comment#1{{}}

\def\qed{{\hfill{\vrule height5pt width3pt depth0pt}\medskip}}

\begin{document}
\begin{center}
{\bf \LARGE  Topological method for symmetric periodic orbits for
maps with a reversing symmetry}

\vskip 0.5cm

\vskip 0.5cm {\large Daniel Wilczak}\footnote{research supported
by Polish State Committee for Scientific
Research grant 2 P03A 041 24 
}  \\
  WSB -- NLU, Faculty of Computer Science,\\
  Department of Computational Mathematics,\\
  Zielona 27, 33-300 Nowy S\c{a}cz, Poland\\
  and \\
  Jagiellonian University, Institute of Computer Science\\
  Nawojki 11, 30-072  Krak\'ow, Poland\\
  e-mail: dwilczak@wsb-nlu.edu.pl \\

\vskip 0.5cm

{\large and}

\vskip 0.5cm

{\large Piotr Zgliczy\'nski}  \\
Jagiellonian University, Institute of Mathematics, \\ Reymonta 4,
30-059 Krak\'ow, Poland \\ e-mail: zgliczyn@im.uj.edu.pl

\end{center}

\vskip 0.5cm

\begin{center}
\today
\end{center}

\begin{abstract}
We present a topological method of obtaining the existence of
infinite number of symmetric periodic orbits for systems with
reversing symmetry. The method is based on covering relations. We
apply the method to a four-dimensional reversible map.
\end{abstract}



\input intro.tex

\input toptool.tex

\input trans.tex

\input rev-sys.tex

\input example.tex

\input numeric.tex

\input ref.tex
\end{document}

%% file: intro.tex
\section{Introduction}
\label{sec:intro}

The goal of this paper is to present a topological method, which
allows to establish, in finite computation, the existence of an
infinite number of symmetric periodic orbits for dynamical systems
with a reversing symmetry.

The role and impact of time reversing symmetries in dynamical
systems has been extensively covered in the literature, see
\cite{la, la1,d1,d2} and references given there.  In this paper we
will restrict ourselves to the discrete time case. Before we
outline our results we need first to introduce a few definitions.

\begin{definition}
\label{def:sym_map} An invertible transformation
$M:\Omega\longrightarrow \Omega$ is called a \emph{symmetry} of
the discrete dynamical system induced by
$f:\Omega\longrightarrow\Omega$ if
\[
M\circ f=f\circ M.
\]
\end{definition}
It is easy to see that a symmetry of a discrete dynamical system
transforms trajectories into trajectories.

Alternatively, a dynamical system may admit transformations
$S:\Omega \rightarrow\Omega$ that map trajectories into other
trajectories reversing the time direction.
\begin{definition}\cite{la1}
An invertible transformation $S:\Omega\longrightarrow\Omega$ is
called a reversing symmetry for the discrete dynamical system
induced by $f:\Omega \longrightarrow\Omega$ if
\[
S\circ f=f^{-1}\circ S,
\]
or equivalently $S^{-1}\circ f\circ S\circ
f=\mathrm{Id}|_{\Omega}$.
\end{definition}

\begin{definition}
Let $\varphi:\mathbb{T}\times\Omega\rightarrow\Omega$ (where $\mathbb{T}%
=\mathbb{R}$ or $\mathbb{T}=\mathbb{Z}$) be a dynamical system.
For each point $x_{0}\in\Omega$, we define its orbit
$\Gamma(x_{0})$ by
\[
\Gamma(x_{0})=\{x \: | \: x=\varphi(t,x_{0})\text{ for some
}t\in\mathbb{T}\}
\]
\end{definition}

To describe the symmetry properties of orbits, following
\cite{gss,la1} we introduce the notion of isotropy subgroup.

\begin{definition}
Consider a dynamical system on a phase space $\Omega$ with a
symmetry (reversing symmetry) $S$.

If $S(\Gamma(x_0))= \Gamma(x_0) $ for some $x_0 \in \Omega$, then
we say that the orbit $\Gamma(x_{0})$ is $S$-symmetric. If the
symmetric orbit $\Gamma(x_0)$ is periodic, then we will say that
$x_0$ is a symmetric periodic point.

If $S$ is clear from the context, then we will often drop it and
speak of a symmetric orbit.
\end{definition}

Let $P:X \to X$ be a map with a reversing symmetry $S$. Let
$\Fix(S)=\{ x \  | \ S(x)=x \}$ be the fixed point set for $S$.

A standard method for studying  symmetric periodic orbits for
systems with  a reversing symmetry is the Fixed Set Iteration
(FSI) method \cite{lq,la1} (also known as DeVogelaere method
\cite{dv}). This method is based on the intersections of the
iterates of $\Fix(S)$ and $\Fix(P\circ S)$ (see \cite[Prop.
1.2.2]{la1} for more details). Below we give its simplified
version.
\begin{theorem}
\label{thm:fsi} If  $y_1 \in \Fix(S)$ and $P^k(y_1) \in \Fix(S)$
for some $k \in \mathbb{Z}_+$, then $y_1$ is a symmetric periodic
point  for $P$ and its  principal period divides $2k$.
\end{theorem}

Theorem~\ref{thm:fsi} shows that to obtain a symmetric periodic
orbit of period $k$, we have to examine  the $k/2$-iterate of $P$,
which is an apparent obstacle for obtaining infinite number of
$S$-symmetric periodic orbits with unbounded periods in a finite
computation.

One way to overcome this problem is the method proposed by Devaney
in \cite{d2}. To describe the Devaney's method, assume that
$P:\mathbb{R}^{2n} \to \mathbb{R}^{2n}$ is a $C^1$-map with an
reversing symmetry $S$ and $\Fix(S)$ is an $n$-dimensional
manifold. Now, if we have $p \in \Fix(S)$, which is a fixed
(periodic) point for $P$ and the unstable manifold of $p$,
$W^u(p)$, intersects transversally $\Fix(S)$, then  by symmetry
$W^s(p)=S(W^u(p))$, hence $W^s(p)$  intersects $\Fix(S)$
transversally at $q$. Now we apply $P^k$ to the neighborhood of
$q$ and it is easy to see, that for $k$ large enough there exists
a disk $D_k$ in $\Fix(S)$, such that $P^k(D_k)$ intersects
$\Fix(S)$ transversally in the vicinity of $p$. Now by
Theorem~\ref{thm:fsi} one obtains periodic points of arbitrary
high periods.

The method proposed in this paper bears some resemblance with the
Devaney's method, as it also is based on  some kind of
transversality, which is preserved under the iteration of the map.
Our method is purely topological and is based on the notion of
covering relation (see \cite{GiZ}) as a tool for the propagation
of the topological transversality. Here we will informally outline
our approach. We say that a cube $N$ $P$-covers a cube $M$,
denoted by $N \cover{P} M$, if the image of the cube $N$ under the
map $P$ is stretched across the cube $M$ in a topologically
nontrivial manner (see Definition \ref{def:covw}). These cubes
together with the corresponding choices of coordinate systems will
be referred to as h-sets (the letter h suggesting the
hyperbolic-like directions). Now in each h-set $N$ we can define
horizontal and vertical disks (see Definitions~\ref{def:hor} and
\ref{def:ver}). In Section~\ref{sec:trans} we prove (see
Theorem~\ref{thm:trans}) that if we have the chain of covering
relations
\begin{equation}
  N=N_0 \cover{P} N_1 \cover{P} \dots \cover{P} N_k=M, \label{eq:chain-cov-rel}
\end{equation}
then for any horizontal disk $H$ in $N_0$ and any vertical disk
$V$ in $N_k$, there exists $x \in H$, such that $P^i(x) \in N_i$
for $i=1,\dots,k$ and $P^k(x) \in V$. Now if $\Fix(S)$ forms a
horizontal disk in $N$ and $\Fix(S)$ forms a vertical disk in $M$,
then from Theorem~\ref{thm:fsi} we obtain symmetric periodic
point. Observe that if we can build a chain of covering relations
(\ref{eq:chain-cov-rel}) linking $N$ and $M$ of arbitrary length,
then we will have symmetric periodic points of arbitrary high
period. For example it was shown that this happens for suitable
2-dimensional Poincar\'e maps for the Michelson system arising
from the Kuramoto-Sivashinsky PDE \cite{W}, the planar restricted
three body problem modelling the motion of the Oterma comet in the
Sun-Jupiter system \cite{WZ2} or the Henon-Heiles Hamiltonian
\cite{az}. In the above mentioned applications we had one unstable
and one stable directions. In this paper we apply our method to a
four-dimensional reversible map to prove the existence of symbolic
dynamics and an infinite number of symmetric periodic orbits. The
main feature, which makes this example interesting,  is the fact
that both stable and unstable directions are two-dimensional. The
proof is computer assisted, i.e., rigorous numerics is was used to
verify assumptions of abstract theorems.

The proposed method was introduced in \cite{W} in the planar case
and direct coverings. In this case the proof of the transversality
theorem (corresponding to Theorem~\ref{thm:trans}) is very simple
and is based on the connectivity argument. Unfortunately this
proof cannot be generalized to higher dimension or to include
coverings induced by inverse mappings.

The content of this paper can be described as follows. In
Section~\ref{sec:topmeth} we define topological notions: h-sets,
covering relations and backcoverings relations. In
Section~\ref{sec:trans} we prove the main transversality theorem
(Theorem~\ref{thm:trans}). In Section~\ref{sec:rev-sys} describe
how it can be applied  to reversible dynamical systems in general.
In Section~\ref{sec:4d-rev} we consider an application of our
method to a  four-dimensional  map with an reversing symmetry. In
Section~\ref{sec:num} we describe how to verify the existence of
covering relations by computer.

%% file: toptool.tex
\section{Topological tools:  h-sets and covering relations}
\label{sec:topmeth}

In this section we present  main topological tools used in this
paper. The crucial notion is that of {\em  covering relation}
\cite{GiZ}.

\subsection{h-sets}

{\bf Notation:} For a given norm in $\mathbb{R}^n$ by $B_n(c,r)$
we will denote an open ball of radius $r$ centered at $c \in
\mathbb{R}^n$. When the dimension $n$ is obvious from the context
we will drop the subscript $n$. Let $S^n(c,r)=\partial
B_{n+1}(c,r)$, by the symbol $S^n$ we will denote $S^n(0,1)$. We
set $\mathbb{R}^0=\{0\}$, $B_0(0,r)=\{0\}$, $\partial
B_0(0,r)=\emptyset$.

For a given set $Z$, by $\inte Z$, $\overline{Z}$, $\partial Z$ we
denote the interior, the closure and the boundary of $Z$,
respectively. For the map $h:[0,1]\times Z \to \mathbb{R}^n$ we
set $h_t=h(t,\cdot)$.  By $\mbox{Id}$ we denote the identity map.
For a map $f$, by $\mbox{dom} (f)$ we will denote the domain of
$f$. Let $f : \Omega \subset {\mathbb R}^n \to {\mathbb R}^n$ be a
continuous map, then we will say that $X \subset \dom(f^{-1})$ if
the map $f^{-1}:X \to {\mathbb R}^n$ is well defined and
continuous. For $N \subset \Omega$, $N$-open and $c \in
\mathbb{R}^n$ by $\deg(f,N,c)$ we denote the local Brouwer degree.
For the properties of this notion we refer the reader to \cite{L}
(see also Appendix in \cite{GiZ}).

\begin{definition} \cite[Definition 1]{GiZ}
\label{def:covrel} A $h$-set, $N$, is the  object consisting of
the following data
\begin{itemize}
 \item $|N|$ - a compact subset of ${\mathbb R}^n$
 \item $u(N),s(N) \in \{0,1,2,\dots\}$, such that $u(N)+s(N)=n$
 \item a homeomorphism $c_N:{\mathbb R}^n \to
   {\mathbb R}^n={\mathbb R}^{u(N)} \times {\mathbb R}^{s(N)}$,
     such that
      \begin{displaymath}
        c_N(|N|)=\overline{B_{u(N)}}(0,1) \times
        \overline{B_{s(N)}}(0,1).
      \end{displaymath}
\end{itemize}
We set
\begin{eqnarray*}
  N_c=\overline{B_{u(N)}}(0,1) \times \overline{B_{s(N)}}(0,1), \\
   N_c^-=\partial \overline{ B_{u(N)}}(0,1) \times
\overline{B_{s(N)}}(0,1) \\
N_c^+=\overline{B_{u(N)}}(0,1) \times
\partial \overline{B_{s(N)}}(0,1) \\
  N^-=c_N^{-1}(N_c^-) , \quad N^+=c_N^{-1}(N_c^+)
\end{eqnarray*}
\end{definition}

Hence a $h$-set, $N$, is a product of two closed balls in some
coordinate system. The numbers, $u(N)$ and $s(N)$, stand for the
dimensions of nominally unstable and stable directions,
respectively. The subscript $c$ refers to the new coordinates
given by homeomorphism $c_N$. Observe that if $u(N)=0$, then
$N^-=\emptyset$ and if $s(N)=0$, then $N^+=\emptyset$. In the
sequel to make notation less cumbersome we will drop the bars in
the symbol $|N|$ and we will use $N$ to denote both the h-sets and
its support.

\begin{definition} \cite[Definition 3]{GiZ}
Let $N$ be a $h$-set. We define a $h$-set $N^T$ as follows
\begin{itemize}
 \item $|N^T|=|N|$
 \item $u(N^T)=s(N)$,  $s(N^T)=u(N)$
 \item We
 define a homeomorphism $c_{N^T}:{\mathbb R}^n \to   {\mathbb R}^n={\mathbb R}^{u(N^T)} \times {\mathbb
R}^{s(N^T)}$,
 by
      \begin{displaymath}
        c_{N^T}(x)= j(c_{N}(x)) ,
      \end{displaymath}
      where $j: {\mathbb R}^{u(N)} \times {\mathbb R}^{s(N)} \to {\mathbb R}^{s(N)} \times {\mathbb R}^{u(N)}$
      is given by $j(p,q)=(q,p)$.
\end{itemize}
\qed
\end{definition}
Observe that $N^{T,+}=N^-$ and $N^{T,-}=N^+$. This operation is
useful in the context of inverse maps.

\subsection{Covering relations}
\begin{definition}\cite[Definition 6]{GiZ}
\label{def:covw} Assume that $N,M$ are $h$-sets, such that
$u(N)=u(M)=u$ and $s(N)=s(M)=s$. Let $f:N \to {\mathbb R}^n$ be a
continuous map. Let $f_c= c_M \circ f \circ c_N^{-1}: N_c \to
{\mathbb R}^u \times {\mathbb R}^s$. Let $w$ be a nonzero integer.
We say that
\begin{displaymath}
  N\cover{f,w} M
\end{displaymath}
($N$ $f$-covers $M$ with degree $w$) iff the following conditions
are satisfied
\begin{description}
\item[1.] there exists a continuous homotopy $h:[0,1]\times N_c \to {\mathbb R}^u \times {\mathbb R}^s$,
   such that the following conditions hold true
   \begin{eqnarray}
      h_0&=&f_c,  \label{eq:hc1} \\
      h([0,1],N_c^-) \cap M_c &=& \emptyset ,  \label{eq:hc2} \\
      h([0,1],N_c) \cap M_c^+ &=& \emptyset .\label{eq:hc3}
   \end{eqnarray}
\item[2.] There exists a  map $A:{\mathbb R}^u \to {\mathbb
R}^u$, such that
   \begin{eqnarray}
    h_1(p,q)&=&(A(p),0), \mbox{ for $p \in \overline{B_u}(0,1)$ and $q \in
    \overline{B_s}(0,1)$,}\label{eq:hc4}\\
      A(\partial B_u(0,1)) &\subset & {\mathbb R}^u \setminus
      \overline{B_u}(0,1).  \label{eq:mapaway}
   \end{eqnarray}
  Moreover, we require that
\begin{displaymath}
  \deg(A,\overline {B_u}(0,1),0)=w,
\end{displaymath}
\end{description}
\end{definition}

Note that in the case $u=0$, if $N \cover{f,w} M$, then $f(|N|)
\subset \inte |M|$ and  $w=1$.

Intuitively,  $N \cover{f} M$ if $f$ stretches  $N$ in the
'nominally unstable' direction, so that its projection onto
'unstable' direction in $M$ covers in topologically nontrivial
manner projection of $M$. In the 'nominally stable' direction $N$
is contracted by $f$. As a result $N$ is mapped across $M$ in the
unstable direction, without touching $M^+$. It is also very
helpful to note  that the degree $w$ in the covering relation
depends only on $A_{|\partial B_u(0,1)}$.

\begin{definition}\cite[Definition 7]{GiZ}
 Assume $N,M$ are $h$-sets, such that $u(N)=u(M)=u$ and
$s(N)=s(M)=s$.  Let $g:\Omega \subset {\mathbb R}^n \to {\mathbb
R}^n$. Assume that $g^{-1}:|M| \to {\mathbb R}^n$ is well defined
and continuous. We say that $N \invcover{g,w} M$ ($N$
$g$-backcovers $M$ with degree $w$) iff $M^T \cover{g^{-1},w}
N^T$.
\end{definition}

The following theorem was proved in \cite{GiZ}.
\begin{theorem} \cite[Theorem 9]{GiZ}
\label{th:top}
 Assume $N_i$, $i=0,\dots,k$, $N_k=N_0$ are
$h$-sets and for each $i=1,\dots,k$ we have either
\begin{equation}
  N_{i-1} \cover{f_i,w_i} N_{i} \label{eq:dirgcov}
\end{equation}
or  $|N_{i}| \subset \dom(f_i^{-1})$ and
\begin{equation}
  N_{i-1} \invcover{f_i,w_i} N_{i}.  \label{eq:invgcov}
\end{equation}

Then there exists a point $x \in \inte |N_0|$, such that
\begin{eqnarray}
   f_i \circ f_{i-1}\circ \cdots \circ f_1(x) &\in& \inte |N_i|, \quad i=1,\dots,k \\
  f_k \circ f_{k-1}\circ \cdots \circ f_1(x) &=& x
\end{eqnarray}
\end{theorem}

Obviously we cannot make any claim about the uniqueness of $x$ in
Theorem~\ref{th:top}.

 Theorem~\ref{th:top} shows that both the direct and the inverse
 covering relation can be treated on the same footing. This
 justifies the following definition.
\begin{definition}
Assume $N,M$ are h-sets and $P$ is a continuous map. We say that
\begin{displaymath}
  N \bicover{P,w} M
\end{displaymath}
If one of the two following conditions is satisfied
\begin{eqnarray*}
  N \subset \dom(P) \quad \mbox{and} \quad N \cover{P,w} M \\
  |M| \subset \dom(P^{-1}) \quad \mbox{and} \quad  N \invcover{P,w} M.
\end{eqnarray*}
\end{definition}
We would like to stress, that the relation $N \bicover{P,w} M$ is
not symmetric.

%% file: trans.tex
\section{The topological transversality theorem }
\label{sec:trans} The goal of this section is to state and prove
the main topological transversality theorem for chain of covering
relations. For this end we need first  to define the notions of
vertical and horizontal disks in an h-set.
\begin{definition}
\label{def:hor} Let $N$ be an $h$-set. Let
$b:\overline{B_{u(N)}}(0,1) \to |N|$ be  continuous and let
$b_c=c_N \circ b$.  We say that $b$
 is \emph{a horizontal disk in $N$} if
there exists a continuous homotopy $h:[0,1]\times
\overline{B_{u(N)}}(0,1) \to N_c$, such that
\begin{eqnarray}
  h_0&=&b_c \label{eq:hor1}\\
  h_1(x)&=&(x,0), \qquad \mbox{for all $x \in \overline{B_{u(N)}}(0,1)$} \\
  h(t,x) &\in& N_c^-, \qquad \mbox{for all $t \in [0,1]$ and $x \in \partial
  \overline{B_{u(N)}}(0,1)$}\label{eq:hor3}
\end{eqnarray}
\end{definition}

\begin{definition}
\label{def:ver} Let $N$ be an $h$-set. Let
$b:\overline{B_{s(N)}}(0,1) \to |N|$ be  continuous and let
$b_c=c_N \circ b$. We say that $b$
 is \emph{a vertical disk in $N$} if
there exists a continuous homotopy $h:[0,1]\times
\overline{B_{s(N)}}(0,1) \to N_c$, such that
\begin{eqnarray*}
  h_0&=&b_c \\
  h_1(x)&=&(0,x), \qquad \mbox{for all $x \in \overline{B_{s(N)}}(0,1)$} \\
  h(t,x) &\in& N_c^+, \qquad \mbox{for all $t \in [0,1]$ and $x \in \partial
  \overline{B_{s(N)}}(0,1)$}.
\end{eqnarray*}
\end{definition}
It is easy to see that  $b$ is the horizontal disk in $N$ iff $b$
is the vertical disk in $N^T$.

\begin{figure}[htpb]
    \centerline{\includegraphics[width=2.5in]{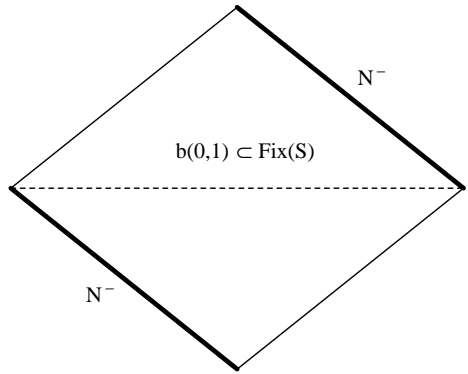}}
    \caption{The curve $b$ is both horizontal and vertical disk in $N$.
      In this example $u(N)=s(N)=1$. \label{fig:symset}}
\end{figure}

We would like to remark here that the horizontal disk in $N$ can
be at the same time also vertical in $N$. An example of such disk
is shown on Fig.~\ref{fig:symset}. In case  homotopies used in the
definitions of horizontal and vertical disks are different. The
existence of such disks, which are both vertical and horizonal
will play very important role in our method for detection of an
infinite number symmetric periodic orbits for maps with reversal
symmetry.

Now we are ready to state and prove the main topological
transversality theorem. A simplified version of this theorem was
given in \cite{W} for the case of one unstable direction and
covering relations chain without backcoverings. The argument in
\cite{W}, which was quite simple and was based on the connectivity
only, cannot be carried over to larger number of unstable
directions or to the situation when both covering and backcovering
relations are present.
\begin{theorem}
\label{thm:trans} Let $k \geq 1$. Assume $N_i$, $i=0,\dots,k$,
are $h$-sets and for each $i=1,\dots,k$ we have either
\begin{equation}
  N_{i-1} \cover{f_i,w_i} N_{i} \label{eq:tran-dirgcov}
\end{equation}
or  $|N_{i}| \subset \dom(f_i^{-1})$ and
\begin{equation}
  N_{i-1} \invcover{f_i,w_i} N_{i}.  \label{eq:tran-invgcov}
\end{equation}
Assume that $b_0$ is a horizontal disk in $N_0$ and $b_e$ is a
vertical disk in $N_k$.

Then there exists a point $x \in \inte |N_0|$, such that
\begin{eqnarray}
   x&=&b_0(t), \quad   \mbox{for some $t \in B_{u(N_0)}(0,1)$} \\
   f_i \circ f_{i-1}\circ \cdots \circ f_1(x) &\in& \inte |N_i|, \quad i=1,\dots,k \\
  f_k \circ f_{k-1}\circ \cdots \circ f_1(x)&=&b_e(z), \quad \mbox{for some $z \in B_{s(N_k)}(0,1)$}
\end{eqnarray}
\end{theorem}
\textbf{Proof:} Without any loss generality we can assume that
\begin{eqnarray*}
  c_{N_i}=\mbox{Id},  \qquad \mbox{for $i=0,\dots,k$}.
\end{eqnarray*}
Then
\begin{eqnarray*}
  f_i=f_{i,c}, \qquad  \mbox{for $i=1,\dots,k$}, \\
  N_i=N_{c,i}, \quad N_i^\pm= N_{i,c}^\pm \quad \mbox{for $i=0,\dots,k$}
 \end{eqnarray*}

We  define $g_i=f^{-1}_i$, for those $i$ for which we have the
back-covering relation $N_{i-1} \invcover{f_i,w_i} N_i$.

Notice that from the definition of covering relation, it follows
immediately that there are $u\geq 0$, $s\geq 0$, such that
$u(N_i)=u$ and $s(N_i)=s$,  for all $i=0,\dots,k$.

The idea of the proof is to rewrite our problem as a zero finding
problem for a suitable map, then to compute its local Brouwer
degree to infer the existence of a solution.

As a tool for keeping track of the occurrences of coverings and
backcoverings, we define the map $\dir:\{1,\dots,k\} \to \{0,1\}$
by $\dir(i)=1$ if $N_{i-1} \cover{f_i,w_i} N_{i}$ and $\dir(i)=0$
if $N_{i-1} \invcover{f_i,w_i} N_{i}$. For $i=1,\dots,k$ let $h_i$
be a homotopy map from the definition of covering relation for
$N_{i-1} \cover{f_i,w_i} N_i$ or $N_{i-1} \invcover{f_i,w_i} N_i$.
In the case of a direct covering (i.e. $\dir(i)=1$), the homotopy
$h_i$ satisfies
  \begin{eqnarray}
      h_i(0,x)&=&f_i(x), \quad \mbox{where $x \in {\mathbb R}^{u+s}$}, \label{eq:cd1} \\
      h_i(1,(p,q))&=&(A_i(p),0), \quad \mbox{where $p \in {\mathbb R}^u$ and $q \in {\mathbb
      R}^s$},
      \label{eq:cd2} \\
      h_i([0,1],N_{i-1}^-) \cap N_i &=& \emptyset, \label{eq:cd3} \\
      h_i([0,1],N_{i-1}) \cap N_i^+ &=& \emptyset.  \label{eq:cd4}
   \end{eqnarray}
In the case of a backcovering (i.e. $\dir(i)=0$), the homotopy
$h_i$ satisfies
\begin{eqnarray}
      h_i(0,x)&=&g_i(x), \quad \mbox{where $x \in {\mathbb R}^{u+s}$}, \label{eq:ci1} \\
      h_i(1,(p,q))&=&(0,A_i(q)), \quad \mbox{where $p \in {\mathbb R}^u$ and $q \in {\mathbb R}^s$}, \label{eq:ci2} \\
      h_i([0,1],N_{i}^+) \cap N_{i-1} &=& \emptyset, \label{eq:ci3} \\
      h_i([0,1],N_{i}) \cap N_{i-1}^- &=& \emptyset. \label{eq:ci4}
\end{eqnarray}
Let $h_t$ and $h_z$ be the homotopies appearing in the definition
of the horizontal and vertical disk for $b_0$ and $b_e$,
respectively.

 It is enough to  prove that there exists $t \in B_u(0,1)$, $z
\in B_s(0,1)$ and $x_i \in \inte N_{i}$ for $i=1,\dots,k$ such
that
\begin{eqnarray}
  b_0(t) &=& x_0, \nonumber \\
  f_i(x_{i-1})&=&x_i , \quad \mbox{if $\dir(i)=1$},  \nonumber \\
  g_i(x_i) &=& x_{i-1}, \quad \mbox{if $\dir(i)=0$},
  \label{eq:peroi} \\
  b_e(z) &=& x_k. \nonumber
\end{eqnarray}

We will treat (\ref{eq:peroi}) as a multidimensional system of
equations to be solved. To this end, let us define
\begin{displaymath}
  \Pi = \overline{B_u}(0,1) \times N_1 \times \dots \times
  N_{k-1} \times \overline{B_s}(0,1)
\end{displaymath}
A point  $x \in \Pi$ will be  represented by
$x=(t,x_1,\dots,x_{k-1},z)$.

We define  a map $F=(F_1,\dots,F_{k}): \Pi \to {\mathbb
R}^{(u+s)k}$ as follows: for $i=2,\dots,k-1$ we set
\begin{displaymath}
  F_i(t,x_1,\dots,x_{k-1},z)=
  \begin{cases}
    x_i - f_i(x_{i-1}) & \text{if $\dir(i)=1$}, \\
    x_{i-1} - g_i(x_i) & \text{if $\dir(i)=0$}.
  \end{cases}
\end{displaymath}
For $i=1$ we set
\begin{eqnarray*}
  F_1(t,x_1,\dots,x_{k-1},z)=
  \begin{cases}
    x_1 - f_1(b_0(t)) & \text{if $\dir(1)=1$}, \\
    b_0(t) - g_1(x_1) & \text{if $\dir(1)=0$}.
  \end{cases}
\end{eqnarray*}
For $i=k$ we define
\begin{eqnarray*}
  F_k(t,x_1,\dots,x_{k-1},z)=
  \begin{cases}
    b_e(z) - f_k(x_{k-1}) & \text{if $\dir(k)=1$}, \\
    x_{k-1} - g_{k}(b_e(z)) & \text{if $\dir(k)=0$}.
  \end{cases}
\end{eqnarray*}

 With this notation, solving  the system (\ref{eq:peroi}) is
equivalent to solving the equation $F(x)=0$ in $\inte \Pi$.

We define a homotopy $H=(H_1,\dots,H_{k}): [0,1] \times \Pi \to
{\mathbb R}^{(u+s)k}$  as follows. For $i=2,\dots,k-1$ we set
\begin{displaymath}
  H_i(\lambda,t,x_1,\dots,x_{k-1},z)=
  \begin{cases}
    x_i - h_i(\lambda, x_{i-1}) & \text{if $\dir(i)=1$}, \\
    x_{i-1} - h_i(\lambda,x_i) & \text{if $\dir(i)=0$}.
  \end{cases}
\end{displaymath}
For $i=1$ we set
\begin{eqnarray*}
  H_1(\lambda,t,x_1,\dots,x_{k-1},z)=
  \begin{cases}
    x_1 - h_1(\lambda,h_t(\lambda,t)) & \text{if $\dir(1)=1$}, \\
    h_t(\lambda,t) - h_1(\lambda,x_1) & \text{if $\dir(1)=0$}.
  \end{cases}
\end{eqnarray*}
For $i=k$ we define
\begin{eqnarray*}
  H_k(\lambda,t,x_1,\dots,x_{k-1},z)=
  \begin{cases}
    h_z(\lambda,z) - h_k(\lambda,x_{k-1}) & \text{if $\dir(k)=1$}, \\
    x_{k-1} - h_{k}(\lambda,h_z(\lambda,z)) & \text{if $\dir(k)=0$}.
  \end{cases}
\end{eqnarray*}

 Notice  that $H(0,x)=F(x)$. The assertion of the theorem is a
consequence of the following two lemmas, which will be proved
after we complete the current proof.

\begin{lemma}
\label{lem:Icont}
  For all $\lambda \in [0,1]$ the local Brouwer degree  $\deg(H(\lambda,\cdot),\inte
  \Pi,0)$ is well defined and does not depend on $\lambda$. Namely,
  for all $\lambda \in [0,1]$ we have
  \begin{displaymath}
    \deg(H(\lambda,\cdot),\inte  \Pi,0)= \deg(H(1,\cdot),\inte  \Pi,0).
  \end{displaymath}
\end{lemma}

\begin{lemma}
\label{lem:Idegcalc}
\begin{displaymath}
  \left| \deg(H(1,\cdot),\inte  \Pi, 0) \right|=|w_1 \cdot w_2 \cdot \dots \cdot w_k |
\end{displaymath}
\end{lemma}

We continue the proof of Theorem~\ref{thm:trans}. Since
$F=H(0,\cdot)$, from the above lemmas it follows immediately that
\begin{displaymath}
    \deg(F,\inte  \Pi,0)= \deg (H(0,\cdot),\inte  \Pi, 0)= \deg(H(1,\cdot),\inte
  \Pi, 0) \neq 0.
\end{displaymath}
Hence there exists $x \in \Pi$ such that $F(x)=0$.
 \qed

\noindent \textbf{Proof of Lemma~\ref{lem:Icont}:}
 From the homotopy
property of the local Brouwer degree (see Appendix in \cite{GiZ})
it is enough to prove that
\begin{equation}
  H(\lambda,x) \neq 0, \quad \mbox{for all $x \in \partial(\Pi)$  and $\lambda \in [0,1]$.}
  \label{eq:iboundok}
\end{equation}

In order to prove (\ref{eq:iboundok}), let us fix
$x=(t,x_1,\dots,x_{k-1},k) \in \partial{\Pi}$. It is easy to see
that one of the following conditions must be satisfied
\begin{eqnarray}
  t &\in& \partial B_u(0,1),  \label{eq:bdhor} \\
  z &\in& \partial B_s(0,1),  \label{eq:bdver} \\
  x_i &\in& N_i^+, \quad \mbox{for some $i=1,\dots,k-1$}  \label{eq:iinentry}, \\
  x_i &\in& N_i^-  \quad \mbox{for some $i=1,\dots,k-1$}  \label{eq:iinexit}.
\end{eqnarray}
We will deal with all above cases separately.

Consider first (\ref{eq:bdhor}). Let us fix $\lambda \in [0,1]$.
Let $x_0=h_t(\lambda,t)$. From Def.~\ref{def:hor}  it follows that
$x_0 \in N_0^-$. There are now two possibilities: either
$\dir(1)=1$ (direct covering) or $\dir(1)=0$ (backcovering).
Assume that $\dir(1)=1$. From condition (\ref{eq:cd3}) it follows
that $h_1(\lambda,x_0) \notin N_1$, hence $H_1(\lambda,x) \neq 0$.
Assume now that $\dir(1)=0$. We have $h_t(\lambda,t) \in N_0^-$
and since by (\ref{eq:ci4})  $h_1(\lambda,N_1) \cap N_0^- =
\emptyset$, hence $H_1(\lambda,x) \neq 0$.

Consider now (\ref{eq:bdver}). Let us fix $\lambda \in [0,1]$ and
let $x_k=h_z(\lambda,z)$. From Def.~\ref{def:ver} it follows that
$x_k \in N_k^+$. Now if $\dir(k)=1$, then from condition
(\ref{eq:cd4}) $h_k(\lambda,N_{k-1}) \cap N_k^+=\emptyset$, hence
$H_k(\lambda,x) \neq 0$. If $\dir(k)=0$, then from (\ref{eq:ci3})
if follows, that $h_k(\lambda,x_k) \notin N_{k-1}$, hence
$H_k(\lambda,x) \neq 0$.

For each of cases (\ref{eq:iinentry}) and (\ref{eq:iinexit}) we
have to consider the following four possibilities
\begin{eqnarray}
  N_{i-1} & \cover{f_i}& N_{i} \cover{f_{i+1}} N_{i+1}, \label{eq:dd} \\
  N_{i-1} & \cover{f_i}& N_{i} \invcover{f_{i+1}} N_{i+1}, \label{eq:di} \\
  N_{i-1} & \invcover{f_i}& N_{i} \cover{f_{i+1}} N_{i+1}, \label{eq:id} \\
  N_{i-1} & \invcover{f_i}& N_{i} \invcover{f_{i+1}} N_{i+1}.  \label{eq:ii}
\end{eqnarray}
Assume first that $x_i \in N_i^+$. If (\ref{eq:dd}) or
(\ref{eq:di}) holds true, then  from (\ref{eq:cd4}) we obtain
\begin{displaymath}
  h_i(\lambda,x_{i-1}) \neq x_i,
\end{displaymath}
 for every $\lambda \in [0,1]$ and every $x_{i-1} \in N_{i-1}$.
 If (\ref{eq:id}) or
(\ref{eq:ii}) is satisfied, then  from (\ref{eq:ci3}) it results
that
\begin{displaymath}
  h_i(\lambda,x_i) \neq x_{i-1},
\end{displaymath}
for  every $\lambda \in [0,1]$ and  every $x_{i-1} \in |N_{i-1}|$.
This proves  if $x_i \in N_i^+$, then $H(\lambda,x) \neq 0$ for
any $\lambda \in [0,1]$.

Assume now that $x_i \in N_i^-$. If (\ref{eq:dd}) or (\ref{eq:id})
holds true, then from (\ref{eq:cd3}) it follows that for every
$\lambda \in [0,1]$ and every $x_{i+1} \in N_{i+1}$ we have
\begin{displaymath}
  h_{i+1}(\lambda,x_{i}) \neq x_{i+1}.
\end{displaymath}
If (\ref{eq:di}) or (\ref{eq:ii}) is satisfied, then from
(\ref{eq:ci4}) we obtain
\begin{displaymath}
  h_{i+1}(\lambda,x_{i+1}) \neq x_i,
\end{displaymath}
for every $\lambda \in [0,1]$ and every $x_{i+1} \in N_{i+1}$.
This proves that if $x_i \in N_i^-$, then $H(\lambda,x) \neq 0$
for any $\lambda \in [0,1]$.
 \qed

\noindent \textbf{Proof of Lemma~\ref{lem:Idegcalc}:}
 Let us represent
$x_i$ for $i=1,\dots,k-1$ as a pair $x_i=(p_i,q_i)$, where $p_i
\in {\mathbb R}^u$ and $q_i \in {\mathbb R}^s$. In this
representation the map
$H(1,t,p_1,q_1,\dots,p_{k-1},q_{k-1},z)=({\tilde p}_1,{\tilde
q}_1,\dots,{\tilde p}_{k},{\tilde q}_{k}) $ has the following form
(for $\alpha=0$)
\begin{itemize}
\item if $i=2,\dots,k-1$  then
\begin{eqnarray}
 \text{if $\dir(i)=1$, then } & &  {\tilde p}_i = (1-\alpha) p_i - A_i(p_{i-1}), \qquad
    {\tilde q}_i = q_i , \\
  \text{if $\dir(i)=0$, then } & & {\tilde p}_i = p_{i-1}, \qquad
  {\tilde q}_i = (1-\alpha)  q_{i-1} - A_i(q_i)
\end{eqnarray}
\item if $i=1$, then
\begin{eqnarray}
\text{if $\dir(1)=1$, then } & &  \tilde{p}_1=(1-\alpha) p_1 -
A_1(t), \qquad \tilde{q}_1=q_1 \\
  \text{if $\dir(1)=0$, then } & &   \tilde{p}_1=t, \qquad \tilde{q}_1=-A_1(q_1)
\end{eqnarray}
\item if $i=k$, then
\begin{eqnarray}
 \text{if $\dir(k)=1$, then } & &  \tilde{p}_k=- A_k(p_{k-1}), \qquad
 \tilde{q}_k=z\\
 \text{if $\dir(k)=0$, then } & &  \tilde{p}_k=p_{k-1}, \qquad \tilde{q}_k=(1-\alpha)q_{k-1} - A_k(z)
\end{eqnarray}
\end{itemize}
The above equations define a homotopy $C:[0,1] \times \Pi \to
{\mathbb R}^{(u+s)k}$. We will show that
$\deg(C(\alpha,\cdot),\inte \Pi,0)$ is independent of $\alpha$ and
then we compute the degree of $C(1,\cdot)$.

\begin{lemma}
For any $\alpha \in [0,1]$
\begin{displaymath}
\deg(C(\alpha,\cdot),\inte \Pi,0) = \deg(C(1,\cdot),\inte \Pi,0).
\end{displaymath}
\end{lemma}
{\bf Proof:} From the homotopy property of the local degree (see
Appendix in \cite{GiZ}), it follows that it is enough to prove
that
\begin{equation}
 C(\alpha,x) \neq 0, \quad \mbox{for all $x \in \partial\Pi $ and $\alpha \in [0,1]$.}
\end{equation}
Let us take $x=(t,p_1,q_1,\dots,p_{k-1},q_{k-1},z) \in
\partial\Pi$. One of the following conditions holds
 true
\begin{eqnarray*}
  t &\in& S^u, \\
  z &\in& S^s, \\
  p_i &\in& S^u, \quad \mbox{for some $i=2,\dots,k-1$} \\
  q_i &\in& S^s, \quad \mbox{for some $i=2,\dots,k-1$}.
\end{eqnarray*}
Assume that $t \in S^u$. If $\dir(1)=1$, then $\|A_1(t)\| > 1$,
hence $\|\tilde{p}_1\| \geq \|A_1(t)\| - \|p_1\| >0$. If
$\dir(1)=0$, then $\tilde{p_1}=t \neq 0$.

Assume that $z \in S^s$. If $\dir(k)=1$, then $\tilde{q}_k=z \neq
0$. If $\dir(k)=0$, then $\| A_k(z)\| > 1$ and we obtain $\|
\tilde{q_k} \| \geq \| A_k(z) \| - \| q_{k-1}\| >0$.

 Assume that $p_i \in S^u$. If $\dir(i+1)=1$, then
$\tilde{p}_{i+1} \neq 0$, because from condition
(\ref{eq:mapaway}) it follows that
\begin{equation}
  \|A_{i+1}(p_{i})\| > 1 \geq \|(1-\alpha) p_{i+1} \|,
\end{equation}
for any $p_{i+1} \in \overline{B_u}(0,1)$.

If $\dir(i+1)=0$, then  obviously $\tilde{p}_{i+1}=p_{i} \neq 0$.

The argument for the case $q_i \in S^s$ is similar.
 \qed

Now we turn to the computation of the degree of $C(1,\cdot)$.
Observe that $C(1,\cdot)$ has the following form: for
$i=2,\dots,k-1$
\begin{eqnarray}
  {\tilde p}_i = - A_i(p_{i-1}),\  {\tilde q}_i = q_i  \quad \mbox{if $\dir(i)=1$} \\
  {\tilde p}_i = p_{i-1},\  {\tilde q}_i = - A_i(q_i)\quad \mbox{if
  $\dir(i)=0$},
\end{eqnarray}
for $i=1$
\begin{eqnarray}
  {\tilde p}_1 = - A_1(t),\  {\tilde q}_1 = q_1  \quad \mbox{if $\dir(1)=1$} \\
  {\tilde p}_1 = t,\  {\tilde q}_1 = - A_1(q_1)\quad \mbox{if
  $\dir(1)=0$},
\end{eqnarray}
and for $i=k$
\begin{eqnarray}
  {\tilde p}_k = - A_{k}(p_{k-1}),\  {\tilde q}_k = z  \quad \mbox{if $\dir(1)=1$} \\
  {\tilde p}_k = t,\  {\tilde q}_k = - A_k(z)\quad \mbox{if
  $\dir(1)=0$}.
\end{eqnarray}

 From the product property of the degree (see Appendix in \cite{GiZ})  it follows that
\begin{eqnarray*}
  |\deg(C(1,\cdot),\Pi,0)| =& \\
    \left| \Pi_{i \in \dir^{-1}(1)}  \deg(-A_i,\overline{B_u}(0,1),0)  \cdot
     \Pi_{i \in \dir^{-1}(0)}
     \deg(-A_i,\overline{B_s}(0,1),0)\right|.
\end{eqnarray*}
In the formula above if $\dir^{-1}(j)=\emptyset$ (for $j=0,1$),
then the corresponding product is set to be equal to $1$.
Similarly  if $u=0$ or $s=0$, then the corresponding product is
also set equal to $1$.

From Collorary 18 in \cite{GiZ} it follows that
\begin{equation}
  \deg(-A,U,0)=(-1)^u \deg(A,U,0).
\end{equation}
This finishes the proof.  \qed

%% file: rev-sys.tex
\section{Reversing symmetry and covering relations.}
\label{sec:rev-sys}

In this section we apply the tools developed in previous sections
to the study of symmetric  periodic orbits for reversible maps.

\begin{theorem}
\label{thm:symcov}Let $S$ be a reversing symmetry for the local
dynamical system induced by the map $P$. Assume that
\begin{displaymath}
M_{0}\bicover{P}M_{1}\bicover{P}M_{2} \cdots\bicover{P}
M_{n-1}\bicover{P}M_{n},
\end{displaymath}
and

\begin{enumerate}
\item  there exists a horizontal disk in $M_0$ contained in
$\Fix(S)$ \item  there exists a vertical disk in $M_n$ contained
in $\Fix(S)$
\end{enumerate}

Then there exists $x\in M_{0}$, such that%
\begin{gather*}
S(x)=x\text{, }P^{2n}(x)=x\\
\text{the orbit of }x\text{ is }S\text{-symmetric}\\
P^{i}(x)\in M_{i}\text{ for }i=1,\dots,n ,\\
P^{n+i}(x)\in S(M_{n-i})\text{ for }i=1,\dots,n.
\end{gather*}
\end{theorem}
\textbf{Proof:} From Theorem~\ref{thm:trans} it follows that there
exists $x \in \Fix(S)$, such that
\begin{eqnarray*}
  P^{i}(x)\in M_{i}\text{ for }i=1,\dots,n ,\\
  P^n(x) \in \Fix(S).
\end{eqnarray*}
The assertion now follows from the reversing symmetry of $P$. \qed

Now we turn our attention to the action of symmetry on h-sets and
covering relations.
\begin{definition}
\label{def:symm}Let $N$ be an h-set in $\mathbb{R}^n$. Let
$L:\mathbb{R}^n \to \mathbb{R}^n$ be a homeomorphism.

We define an h-set $L*N$ as follows
\begin{itemize}
\item $|L*N|=L(|N|)$,
\item $u(L*N)=u(N)$ and $s(L*N)=s(N)$,
\item $c_{L*N}=c_N \circ L^{-1}$.
\end{itemize}

We define an h-set $L^T*N$ by
\begin{displaymath}
  L^T*N= (L*N)^T
\end{displaymath}
\end{definition}
Informally speaking, $L*N$ is just a natural symmetric image of
$N$ and $L^T*N$ is the symmetric image of $N$, but we additionally
switch the 'expanding' and 'contracting' directions.

We have the following
\begin{lemma}
\label{lem:sym-cov} Let $S:\mathbb{R}^{2n} \to \mathbb{R}^{2n}$ be
a reversing symmetry for a map $P$ and $N \bicover{P} M$, then
$S^T*M \bicover{P} S^T * N$.
\end{lemma}
\textbf{Proof:} From the definition of covering relations and
reversing symmetry it follows immediately, that
\begin{eqnarray*}
  \mbox{if} \quad N \cover{P} M, \quad \mbox{then} \quad S^T*M
  \invcover{P} S^T*N \\
   \mbox{if} \quad N \invcover{P} M, \quad \mbox{then} \quad S^T*M
  \cover{P} S^T*N
\end{eqnarray*}
\qed

The reversing symmetry maps, in a natural way, horizontal disks to
vertical disks and vice versa. Namely, we have the following
obvious lemma.
\begin{lemma}
\label{lem:sym-slice} Let $S:\mathbb{R}^{n} \to \mathbb{R}^{n}$  a
homeomorphism, $N$ be an h-set and $\gamma$ be a horizontal
(vertical) disk in $N$.

Then $S(\gamma)$ is a vertical (horizontal) disk in $S^T*N$.
\end{lemma}

In the context of proving  the existence of an infinite number of
symmetric periodic orbits  symmetric h-sets are of special
importance.
\begin{definition}
Let $S$ be a reversing symmetry. We say that an h-set $N$ is
\emph{$S^T$-symmetric} if $S^T*N=N$.
\end{definition}
It is easy to see, that if $N$ is $S^T$-symmetric h-set, then
$u(N)=s(N)$ and  the dimension of the phase space $n=u(N)+s(N)$
must be even.

 The following theorem, which is an easy consequence
of Theorem~\ref{thm:symcov}, illustrates  our method of  proving
of the existence of an infinite number of symmetric periodic
orbits.
\begin{theorem}
\label{thm:wil-method}Let $S$ be a reversing symmetry for the
local dynamical system induced by the map $P$. Let $M_i$ for
$i=0,1,2,3$ be h-sets and $M_i \cap M_j = \emptyset$ for $i \neq
j$. Assume that
\begin{gather*}
M_{0}\bicover{P}M_{1}\bicover{P}M_{2}, \\
M_0 \bicover{P} M_0 \\
\text{$M_0$ and $M_2$ are $S^T$-symmetric } \\
  M_3 = S^T*M_1\\
   \text{ there exists a horizontal disk in $M_0$ contained in
$\Fix(S)$}  \\
\text{ there exists a horizontal disk in $M_2$ contained in
$\Fix(S)$}
\end{gather*}

Then for any sequence $(\alpha_0,\alpha_1,\dots,\alpha_n) \in
\{0,1,2,3\}^{n+1}$, satisfying the following conditions
\begin{gather*}
   \alpha_0 \in \{0,2\}, \quad \alpha_n \in \{0,2\} \\
   \text{if } \alpha_i=0, \text{ then  } \quad \alpha_{i+1} \in \{0,1\}   \\
   \text{if } \alpha_i=1, \text{ then  } \quad  \alpha_{i+1}=2   \\
   \text{if } \alpha_i=2, \text{ then  } \quad  \alpha_{i+1}=3   \\
   \text{if } \alpha_i=3, \text{ then  } \quad  \alpha_{i+1}=1,
\end{gather*}
there exists $x\in |M_{\alpha_0}|$, such that%
\begin{gather*}
S(x)=x\text{, }P^{2n}(x)=x\\
\text{the orbit of }x\text{ is }S\text{-symmetric}\\
P^{i}(x)\in |M_{i}|\text{ for }i=1,\dots,n ,\\
P^{n+i}(x)\in S(|M_{n-i}|)\text{ for }i=1,\dots,n.
\end{gather*}
\end{theorem}
\textbf{Proof:} From Lemma~\ref{lem:sym-slice} it follows that in
$M_0$ and $M_2$ there exist vertical disks contained in $\Fix(S)$.
We can now apply Theorem~\ref{thm:symcov} to the chain of covering
relations
\begin{displaymath}
  M_{\alpha_0} \bicover{P} M_{\alpha_1} \bicover{P} \cdots
  \bicover{P} M_{\alpha_n}
\end{displaymath}
 \qed

%% file: example.tex
\section{One four-dimensional reversible example.}

\label{sec:4d-rev} In this section we present the application of
the method introduced throughout the paper to a four-dimensional
reversible map. As a consequence we obtain the existence of
chaotic dynamics and the existence of an infinite number of
symmetric periodic orbits for a certain iteration of such map. The
proof is computer assisted, i.e., rigorous numerics is used to
verify assumptions of abstract theorems. The main feature, which
makes this example interesting  is the fact that both stable and
unstable directions are two-dimensional and the map itself is not
close to a product of two two-dimensional maps, with unstable
dimension each.

\subsection{An example of four-dimensional reversible map.}
Let $f:\mathbb R^n\to\mathbb R^n$, $n>0$ be a continuous map,
$F:\mathbb R^{2n}\to\mathbb R^{2n}$ be a map defined by
\begin{equation*}
F\left[\begin{matrix}x\\y\end{matrix}\right] =
\left[\begin{matrix}-y+f(x)\\x\end{matrix}\right]
\end{equation*}
and let $S(x,y)=(y,x)$. It is straightforward to show that $S\circ
F\circ S \circ F=\mathrm{Id}$. Therefore $F$ is a reversible
homeomorphism of $\mathbb R^{2n}$. In suitable coordinates we may
rewrite $F$ as follows
\begin{equation}\label{eq:revF}
F\left[\begin{matrix}x\\y\end{matrix}\right] =
\left[\begin{matrix}-y+\frac{1}{2}f(x+y)\\x+\frac{1}{2}f(x+y)\end{matrix}\right]
\end{equation}
and the reversing symmetry $S$ will be given by $S(x,y)=(-x,y)$.

Let us fix $n=2$ and let $f:{\mathbb R}^2\to {\mathbb R}^2$ be
defined by
\begin{equation}\label{eq:f}
    f\left[\begin{matrix}x_1\\x_2\end{matrix}\right] =
    \left[\begin{matrix}x_1(1-x_1)+4-x_2\\x_2(1-x_2)+4+x_1\end{matrix}\right].
\end{equation}
In the remainder of this section we will investigate the map $F$
given by (\ref{eq:revF}), where $f$ is as above. The map
$S(x_1,x_2,y_1,y_2)=(-x_1,-x_2,y_1,y_2)$ is a reversing symmetry
of $F$.

It is easy to verify that each solution of
\begin{equation}\label{eq:solvefix}
F(x_1,x_2,y_1,y_2)=(x_1,x_2,y_1,y_2)
\end{equation}
satisfies $x_1=x_2=0$, hence all fixed points of $F$ are
symmetric. Solving of Eq.~(\ref{eq:solvefix}) leads to the
following system of equations
\begin{equation*}
\begin{cases}
    y_1^2+(y_2+1)^2=9\\
    (y_1+1)^2-y_2^2=1
\end{cases}
\end{equation*}
\begin{figure}[htbp]
    \centerline{\includegraphics[width=2.2in]{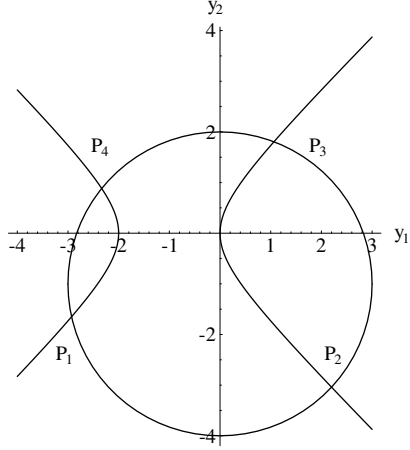}}
    \caption{The location of the fixed points of $F$.\label{fig:fixpoints}}
\end{figure}
describing an intersection of a circle with a hyperbola in four
points as was shown in Fig.~\ref{fig:fixpoints}. Fixed points
$P_1, P_2$ are hyperbolic with real eigenvalues, the fixed point
$P_3$  is hyperbolic with four complex eigenvalues, the point
$P_4$  possesses two complex eigenvalues on the unit circle and
two real eigenvalues, i.e., it is of elliptic-hyperbolic type.

 The fixed points may be exactly computed (for example using
Mathematica). However, for our further consideration it is
sufficient to use two approximate fixed points, which we will
still denote by $P_i$, given by
\begin{eqnarray}
    P_1&=&(0,0,-2.9288690017630725,-1.649404627725545), \label{eq:P1} \\
    P_2&=&(0,0, 2.199939462565084, -3.0396731015162355). \label{eq:P2}
\end{eqnarray}
We will show that in the vicinity of $P_1$ and $P_2$ the map $F^7$
has symbolic dynamics on two  symbols  and there exist an infinite
number of symmetric periodic points of an arbitrary large period.

Before we proceed with the statement of main results for $F$ we
need to discuss how we  represent h-sets.

\subsection{Representation
of h-sets in $\mathbb R^n$.} To define a h-set $N$ we need to
specify a homeomorphism $c_N$ of $\mathbb R^n$ and two numbers
$u(N)$ and $s(N)$ (see Def.~\ref{def:covrel}). Since we will use
the computer in order to verify covering relations, the
homeomorphism must be representable by the machine. The simplest
case is to take an affine map.  We will use the maximum norm on
${\mathbb R}^n$, i.e. $\|x \|=\max_i |x_i|$ and we treat vectors
as columns.

Let $x,u_1,\ldots,u_k,s_1,\ldots,s_{n-k}\in {\mathbb R}^n$, $0\leq
k\leq n$ be such that the vectors
$u_1,\ldots,u_k,s_1,\ldots,s_{n-k}$ are linearly independent. We
define a matrix $M \in {\mathbb R}^{n \times n}$ by
\begin{equation*}
M=[u_1,\ldots,u_k,s_1,\ldots,s_{n-k}].
\end{equation*}
We define a h-set
\begin{equation*}
N=\mathfrak{h}(x,u_1,\ldots,u_k,s_1,\ldots,s_{n-k})
\end{equation*}
as follows
\begin{eqnarray*}
    u(N) &=& k,\\
    s(N) &=& n-k,\\
    |N| &=& M(\overline{B_n}(0,1))+x =
    M(\overline{B_{u(N)}}(0,1)\times\overline{B_{s(N)}}(0,1))+x,\\
    c_N(v) &=& M^{-1}(v-x),\quad \text{for }v\in\mathbb R^n
\end{eqnarray*}
Hence, the h-set $N$ defined above is a parallelepiped centered at
$x$.

In the sequel we will work in $\mathbb{R}^{2k}$ with
$s(N)=u(N)=k$. In this case we will use also the notation
$N=\mathfrak{h}(x,M)$, where  $M$ is a linear  isomorphism of
$\mathbb R^{2k}$. The first $k$ columns of $M$ correspond to
unstable directions and the last $k$ columns of $M$ correspond to
stable directions.

\subsection{Important h-sets and covering relations between them}
\begin{figure}[htbp]
\centerline{
\includegraphics[width=2.5in]{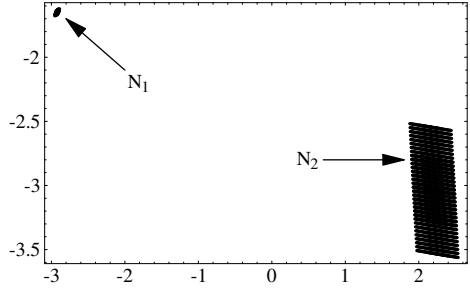}}
\caption{The projection of the sets $|N_1|$ and $|N_2|$ onto
$(y_1,y_2)$ coordinates.\label{fig:n1n2}}
\end{figure}

As was mentioned before $P_1$ and $P_2$ are good numerical
approximations to two hyperbolic fixed points with two-dimensional
stable and unstable manifolds.  We choose $u_j^i,s_j^i \in
\mathbb{R}^4$, $i,j=1,2$ to be a good numerical approximations of
unstable and stable eigenvectors of $DF(P_i)$. Put
\begin{equation}\label{eq:eigenvectors}
\begin{split}
    u_1^1=\left[
    \begin{matrix}
        0.527847408170044\\
        0.254065286036574\\
        0.730261232439584\\
        0.351491787265563
     \end{matrix}\right],\quad
    u_1^2=\left[
    \begin{matrix}
        -0.05726452423754\\
        0.594572575636284\\
        -0.0768865282444865\\
        0.7983061369797889
    \end{matrix}\right],\\
    u_2^1=\left[
    \begin{matrix}
        0.233876807615845\\
        0.485903716548415\\
        0.365235930520818\\
        0.758816138574061
    \end{matrix}\right],\quad
    u_2^2=\left[
    \begin{matrix}
        0.8918103319483236\\
        -0.0858921121857865\\
        0.4421352370808943\\
        -0.0425829663821858)
        \end{matrix}\right]
\end{split}
\end{equation}
and put $s_i^j=S(u_i^j)$, $i,j=1,2$. We define matrices $M_i$ for
$i=1,2$ by
\begin{eqnarray}
    M_1 = 0.012[u_1^1,u_2^1,s_1^1,s_2^1],\quad
    M_2 = 0.31[u_1^2,u_2^2,s_1^2,s_2^2].  \label{eq:defMi}
\end{eqnarray}
We define two $S^T$ symmetric h-sets centered at $P_1$ and $P_2$
by
\begin{eqnarray}\label{eq:defN1N2}
    N_1=\mathfrak h(P_1,M_1),\quad N_2=\mathfrak
    h(P_2,M_2).
\end{eqnarray}
The important remark is that the sets $|N_1|$ and $|N_2|$ are
disjoint. The projection of $|N_1|$ and $|N_2|$ onto $(y_1,y_2)$
coordinates is presented in Fig.~\ref{fig:n1n2}. The following
lemma was proved  with a computer assistance.
\begin{lemma}\label{lem:symcover}
The following covering relations hold: $$N_1\cover{F,1}N_1,\quad
N_2\cover{F,-1}N_2.$$
\end{lemma}
The details of the proof will be presented in
Section~\ref{sec:num}.
\begin{rem}
From above lemma and Theorem 2 it follows immediately, that there
exists fixed points in $N_1$ and $N_2$. However, we cannot claim
that they are unique in $N_i$ nor $S$-symmetric.
\end{rem}

Now we will construct a dynamical link between sets $N_1$ and
$N_2$. More precisely, we construct a chain of covering relations
connecting $N_1$ and $N_2$. For this purpose we look for a point
$Q_1$ in the neighborhood of $P_1$ on the unstable manifold of
$P_1$ and such that $F^k(Q_1)$ (here $k=10$) is close to $P_2$.
Then we define additional points $Q_i$ by taking some forward
iterates of $Q_1$.

To be specific we set
\begin{eqnarray*}
    Q_1 &= &P_1+0.0330092432u_1^1 - 0.048949u_1^2 +0.0004931s_1^1,\\
    Q_2 &= &F^4(Q_1),\\
    Q_3 &= & F(Q_2).
\end{eqnarray*}
We have
\begin{eqnarray*}
    F^{-1}(Q_1) \in |N_1|, \\
    \|F^{-1}(Q_1) -P_1\|<0.006,\\
    \|F^{10}(Q_1)-P_2\|<0.001.
\end{eqnarray*}
The points $Q_1$, $Q_2$ and $Q_3$ will be used as the centers of
new h-sets. We define
\begin{eqnarray*}
    H_1&=&\mathfrak h(Q_1,.001u_1^1,.00175u_2^1,.005s_1^1,0.005s_2^1),\\
    H_2&=&\mathfrak h(Q_2,.28u_1^2,.28u_2^2,.2s_1^2,.38s_2^2),\\
    H_3&=&\mathfrak h(Q_3,.15u_1^2,.15u_2^2,.12s_1^2,.42s_2^2).
\end{eqnarray*}
\begin{lemma}\label{lem:covchain}
The following covering relations hold
\begin{equation*}
N_1\cover{F,1}H_1\cover{F^4,-1}H_2\cover{F,-1}H_3\cover{F,-1}N_2.
\end{equation*}
\end{lemma}
The details of the proof will be presented in
Section~\ref{sec:num}.

Let us comment briefly about the spatial relations between $N_i$
and $H_j$.  The set $H_1$ is close to $N_1$, this is the reason
for using the same stable and unstable directions in the
definitions of these sets. Sets $H_2$ and $H_3$ are close to $N_2$
and as in the previous case we used the same stable and unstable
directions for them.

\begin{figure}[htbp]
\centerline{\includegraphics[width=2.2in]{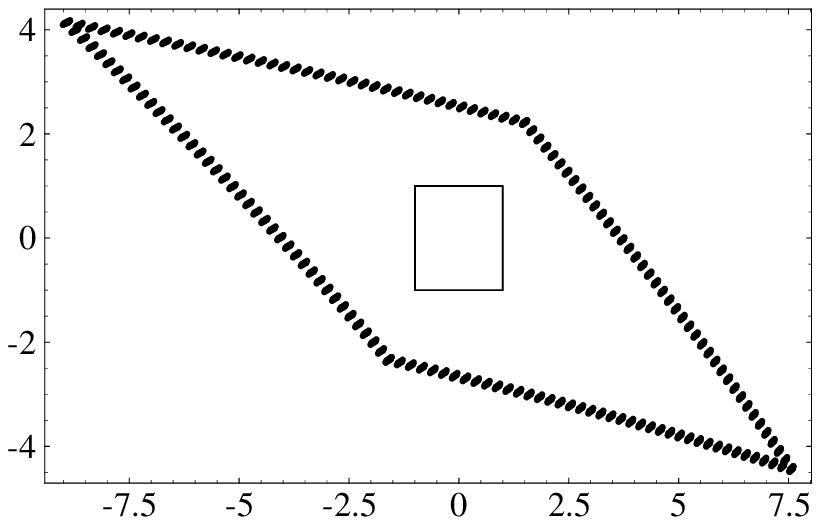}\quad\includegraphics[width=2.2in]{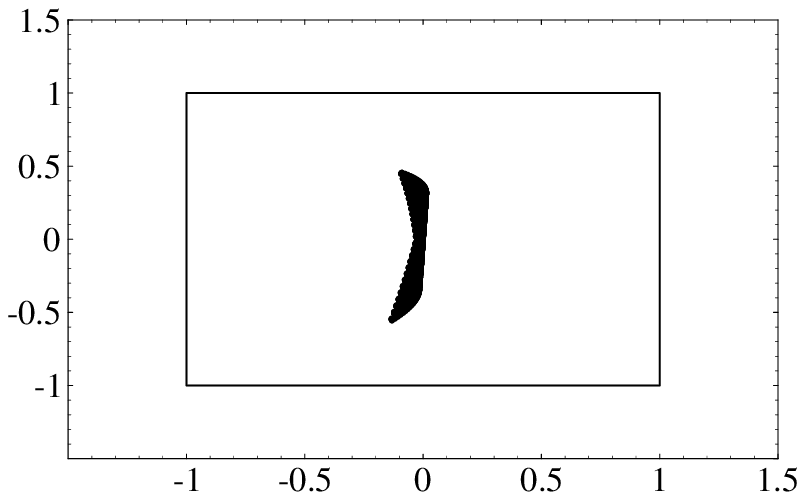}}
\caption{The numerical evidence of the existence of covering
relation $H_1\cover{F^4,-1}H_2$,  (left) the set
$F_c^4((H_1)_c^-)$ projected onto unstable directions of $(H_2)_c$
is outside the unit ball, (right) the set $F_c^4((H_1)_c)$
projected onto stable directions of $(H_2)_c$ is inside the unit
ball.\label{fig:h1image}}
\end{figure}

The numerical evidence of the existence of covering relation
$H_1\cover{F^4,-1}H_2$ is presented in Fig.~\ref{fig:h1image}.

\subsection{Chaotic dynamics of $F$.}
\begin{theorem}
The discrete dynamical system induced by the map $F^7$ is
semiconjugated with the full shift on two symbols, i.e. for an
arbitrary $(i_j)_{j\in\mathbb Z}\in\{1,2\}^{\mathbb Z}$ there
exists a point $x_0\in|N_{i_0}|$ such that
\begin{equation}\label{eq:traj}
F^{7j}(x_0)\in|N_{i_j}|,\quad j\in\mathbb Z.
\end{equation}
Moreover, if the sequence $(i_j)_{j\in\mathbb Z}$ is periodic,
then the point may be chosen as a periodic point for $F^7$ with
the same principal period.
\end{theorem}
\textbf{Proof:} From Lemma~\ref{lem:covchain} we obtain that
\begin{equation}\label{cov:N1N2}
N_1\cover{F,1}H_1\cover{F^4,-1}H_2\cover{F,-1}H_3\cover{F,-1}N_2.
\end{equation}
Since the h-sets $N_1$ and $N_2$ are symmetric (by their
definition), the reversing symmetry property of $F$ implies
\begin{gather}\label{cov:N2N1}
\begin{split}
    N_2=S^T\star N_2
    \invcover{F,-1}S^T\star H_3
    \invcover{F,-1}S^T\star H_2\\
    S^T\star H_2\invcover{F^4,-1}S^T\star H_1
    \invcover{F,1}S^T\star N_1 = N_1.
\end{split}
\end{gather}
From Lemma~\ref{lem:symcover} we get
\begin{eqnarray}
N_1\overbrace{\cover{F,1}N_1\cover{F,1}\cdots\cover{F,1}N_1}^{\text{$7$
times}},\\
N_2\overbrace{\cover{F,-1}N_2\cover{F,-1}\cdots\cover{F,-1}N_2}^{\text{$7$
times}} .\label{cov:N2N2}
\end{eqnarray}
Let $(i_j)_{j\in\mathbb Z}\in\{1,2\}^{\mathbb Z}$ be a periodic
sequence od symbols, i.e. $i_{j+k}=i_{j}$ for $j\in\mathbb Z$ and
a certain $k>0$. Let
\begin{equation}\label{eq:perchain}
N_{i_0}\bicover{}N_{i_1}\cdots\bicover{}N_{i_k}=N_{i_0}
\end{equation}
be a periodic sequence of covering relations, where by
$N_{i_j}\bicover{}N_{i_{j+1}}$ we mean a corresponding sequence
(\ref{cov:N1N2}-\ref{cov:N2N2}). Now, Theorem~\ref{th:top} implies
that there exists a $k$-periodic point $x_0\in|N_{i_0}|$ for
$F^7$, such that assertion (\ref{eq:traj}) is satisfied.

Let $(i_j)_{j\in\mathbb Z}\in\{1,2\}^{\mathbb Z}$ be a nonperiodic
sequence of symbols. We define the periodic sequences
\begin{equation*}
(i_j^k)_{j\in\mathbb
Z}=\{\ldots,\{i_{-k},\ldots,i_0,\ldots,i_k\},\{i_{-k},\ldots,i_0,\ldots,i_k\},\ldots\}\in\{1,2\}^\mathbb
Z
\end{equation*}
where $i_{0}^k=i_0$, $k>0$. Now, for $k>0$ we can find
$(2k+1)-$periodic point $x_k$ for $F^7$ such that
$F^{7j}(x_k)\in|N_{i_j^k}|$, $k>0$, $j\in\mathbb Z$. Since
$N_{i_0}$ is a compact set, we can find a condensation point
$x_0\in|N_{i_0}|$ of $\{x_k\}_{k>0}$. Obviously, $x_0$ satisfies
assertion (\ref{eq:traj}).
 \qed

\subsection{Symmetric periodic points for $F$.}

\begin{theorem}\label{thm:exsym}
There exists an infinite number of symmetric periodic points for
$F$ with an arbitrary large principal periods.
\end{theorem}
The proof of Theorem~\ref{thm:exsym} is a direct consequence of
the following lemma.

\begin{lemma}\label{lem:exsym}
Let
\begin{equation*}
V_0\bicover{g_0}V_1\bicover{g_1}\cdots\bicover{g_{k-1}}V_k
\end{equation*}
be a sequence of covering relations, where
\begin{eqnarray*}
    V_0,V_k\in\{N_1,N_2\},\\
    V_i\in\{N_1,N_2,H_1,H_2,H_3\},\quad \text{for }
    j=1,\ldots,k-1,\\
    g_i=F \text{ or } g_i=F^4,\quad \text{for }i=0,\ldots,k-1.
\end{eqnarray*}
Then there exists a symmetric periodic point $x_0$ of $F$, such
that
\begin{eqnarray*}
    x_0\in|V_0|\cap\Fix(S),\\
    (g_i\circ\cdots\circ g_1\circ g_0)(x_0)\in|V_{i+1}|,\quad
    i=0,\ldots,k-2,\\
    (g_{k-1}\circ\cdots\circ g_1\circ
    g_0)(x_0)\in|V_k|\cap\Fix(S)
\end{eqnarray*}
\end{lemma}
\textbf{Proof:}
Recall, that the sets $N_1$, $N_2$ are symmetric and defined by
vectors $u_i^j,s_i^j$, $i,j=1,2$ -- see(\ref{eq:eigenvectors}) and
(\ref{eq:defN1N2}). For $i=1,2$ we define a map
$b_i:\overline{B_2}(0,1)\longrightarrow |N_i|$ by
\begin{eqnarray*}
    b_i(p,q) & = & M_i [p,q,p,q]^T +P_i, \quad i=1,2.
\end{eqnarray*}
We will show, that $b_i$ is a horizontal disk in $N_i$.

First observe  that
\begin{equation}
  b_{i,c}=c_{N_i} \circ b_i (p,q) = (p,q,p,q).
\end{equation}
We define the homotopy $h_i$ by
\begin{equation*}
    h_i(t,p,q) = (p,q,(1-t)p,(1-t)q).
\end{equation*}
It is easy to see, that conditions  (\ref{eq:hor1}-\ref{eq:hor3})
from Definition~\ref{def:hor} are satisfied. Namely, we have
\begin{eqnarray*}
    h_i(0,p,q)  &= &(p,q,p,q) = (b_i)_c(p,q), \quad \text{for } (p,q)\in\overline{B_2}(0,1),\\
    h_i(1,p,q)  &= &(p,q,0,0), \quad \text{for } (p,q)\in\overline{B_2}(0,1),\\
    h_i(t,p,q) &\in  &(N_i)_c^-=\bd\overline{B_2}(0,1)\times\overline{B_2}(0,1), \quad\text{for }
    (p,q)\in\bd\overline{B_2}(0,1), t\in[0,1].
\end{eqnarray*}
This proves that $b_i$ is a horizontal disk in $N_i$.

Let us remind the reader that $S(P_i)=P_i$ and $S(u^i_j)=s^i_j$.
Hence we obtain
\begin{eqnarray*}
  S(b_i(p,q))=S(k_i(p u^i_1 + q u^i_2 + p s^i_1 + q s^i_2) + P_i)= \\
    =k_i\left(p S(u^i_1) + q S(u^i_2) + p S(s^i_1) + q S(s^i_2)\right) +
    S(P_i)=b_i(p,q).
\end{eqnarray*}
where $k_1=0.012$, $k_2=0.31$ are the coefficients used in
(\ref{eq:defMi}) to define $N_i$. This proves that $b_i$ are
contained in $\Fix(S)$.

We will show that $b_i$ is also vertical disk in $N_i$. Namely, if
follows from  Lemma~\ref{lem:sym-slice}  that $S(b_i)=b_i$ is a
vertical disk in $S^T*N_i=N_i$.

Since $V_0,V_k\in\{N_1,N_2\}$ we conclude that there exists a
horizontal disk contained in $\Fix(S)\cap|V_0|$ and there exists
vertical disk contained in $\Fix(S)\cap|V_k|$. Now, the assertion
is a direct consequence of Theorem~\ref{thm:trans}. \qed

\noindent\textbf{Proof of Theorem~\ref{thm:exsym}:} The assertion
follows from the fact that the sets $|N_1|$ and $|N_2|$ are
disjoint and we can construct an arbitrary number of sequences
satisfying assumptions of Lemma~\ref{lem:exsym}, for example
\begin{equation*}
    N_1\overbrace{\cover{F,1}N_1  \cover{F,1} \dots \cover{F,1} N_1 }^{\text{$k$ \rm times}}
    \cover{F,1}H_1\cover{F^4,-1}H_2\cover{F,-1}H_3\cover{F,-1}H_4\cover{F,-1}N_2,
\end{equation*}
for $k>0$. \qed

%% file: numeric.tex
\section{How to verify covering relations with a computer
assistance.}\label{sec:num}

In this section we discuss some numerical aspects of the
verification of covering relations.

Let $N,M$ be a h-sets in $\mathbb R^n$ such that $u(N)=u(M)=u$ and
$s(N)=s(M)=s$ and let $f:|N|\to\mathbb R^n$ be a continuous. In
order to prove that the covering relation $N\cover{f}M$ holds, it
is necessary to find the homotopy $h:[0,1]\times N_c\to\mathbb
R^u\times \mathbb R^s$ and a map $A:\mathbb R^u\to\mathbb R^u$
satisfying conditions (\ref{eq:hc1}-\ref{eq:hc4}).

Since in our example $f_c=c_M\circ f\circ c_N^{-1}$ is a
diffeomorphism we can try to find a homotopy between $f_c$ and its
derivative computed in the center of set and  projected onto
unstable directions, i.e., we define
\begin{equation*}
A:\mathbb R^u\ni p \to \pi_u(Df_c(0)(p,0))\in\mathbb R^u,
\end{equation*}
where $\pi_u:\mathbb R^n\to\mathbb R^u$ is a projection onto first
$u$ variables. We require  $A$ to be an isomorphism, then we have
\begin{equation*}
\deg(A,\overline{B_u}(0,1),0)=\sgn(\det A)=\pm1.
\end{equation*}
Now we define the homotopy between $f_c$ and $(p,q)\to(A(p),0)$ by
\begin{equation}\label{eq:convhom}
h(t,p,q) = (1-t)f_c+t(A(p),0),\quad \text{for }(p,q)\in
\overline{B_u}(0,1)\times\overline{B_s}(0,1).
\end{equation}
Obviously the homotopy (\ref{eq:convhom}) satisfies conditions
(\ref{eq:hc1}) and (\ref{eq:hc4}). We need to check  if the
homotopy (\ref{eq:convhom})  satisfies conditions
(\ref{eq:hc2}-\ref{eq:hc3}).

Below we describe precise algorithms.

\begin{definition}
Let $U\subset \mathbb R^n$ be a bounded set. We say that $\mathcal
G\subset 2^{\mathbb R^n}$ is a grid of $U$ if
\begin{enumerate}
\item $\mathcal G$ is a finite set
\item $U\subset \bigcup_{G\in\mathcal G}G$
\item each $G\in\mathcal G$ can be represent in a computer
\end{enumerate}
\end{definition}
\begin{definition}
Let $U\subset\mathbb R$ be a bounded set. By $(U)_I$ we denote the
interval enclosure of the set, i.e., the set $(U)_I$ is the
smallest representable interval containing $U$ or
$[-\infty,\infty]$ if there is not a representable interval
containing $U$.

Let $U\subset \mathbb R^n$ be a bounded set. By $(U)_I$ we denote
$(\pi_1(U))_I\times\cdots\times (\pi_n(U))_I$ where $\pi_i$ is a
projection onto $i$-th variable.
\end{definition}
In the algorithms presented bellow all computations are performed
in interval arithmetic \cite{Mo}.

First we discuss how we check condition (\ref{eq:hc2}).
\begin{alg}\label{alg:unstable}
\end{alg}
    \textbf{function} ComputeUnstableWall($\mathcal G_1$ : grid,
    $\mathcal G_2$ : grid) : \textbf{bool}\\
    \textbf{var}\\
    \mbox{}\qquad $X$, $LX$, $Z$ : representable sets;\\
    \textbf{begin}\\
        \mbox{}\qquad\textbf{foreach} $G_1\in\mathcal G_1$\\
        \mbox{}\qquad\textbf{begin}\\
            \mbox{}\qquad\qquad $X:=G_1\times 0$; // $0\in\mathbb R^s$\\
            \mbox{}\qquad\qquad $LX: = (Df_c(0)(X))_I$;\\
            \mbox{}\qquad\qquad\textbf{foreach} $G_2\in\mathcal G_2$\\
            \mbox{}\qquad\qquad\textbf{begin}\\
                \mbox{}\qquad\qquad\qquad $Z:=(f_c(G_1\times G_2)\cup LX)_I$;\\
                \mbox{}\qquad\qquad\qquad \textbf{if not} $\pi_u(Z)\subset \mathbb R^u\setminus \overline{B_u}(0,1)$ \textbf{return False};\\
            \mbox{}\qquad\qquad\textbf{end;}\\
        \mbox{}\qquad\textbf{end;}\\
        \mbox{}\qquad\textbf{return True;}\\
    \textbf{end.}
\begin{lemma}
Assume $N$, $M$ be a h-sets and $f:|N|\to\mathbb R^n$ be such that
$f_c$ is smooth. Let $\mathcal G_1$ be a grid of $\bd B_u(0,1)$
and let $\mathcal G_2$ be a grid of $\overline{B_s}(0,1)$. If
Algorithm~\ref{alg:unstable} is called with arguments $(\mathcal
G_1, \mathcal G_2)$ and returns \textbf{\em True} then the
homotopy defined in (\ref{eq:convhom}) satisfies condition
(\ref{eq:hc2}).
\end{lemma}
\textbf{Proof:} Let $(p,q)\in N_c^-$. Since $\mathcal G_1$ is a
grid of $\bd B_u(0,1)$ and $\mathcal G_2$ is a grid of
$\overline{B_s}(0,1)$ then
\begin{equation*}
\mathcal G_1\times \mathcal G_2:=\{G_1\times G_2\ |\
G_1\in\mathcal G_1, G_2\in\mathcal G_2\}
\end{equation*}
is a grid of $N_c^-$. Therefore $(p,q)\in G_1\times G_2$ for some
$G_1\in\mathcal G_1$, $G_2\in\mathcal G_2$. Since the Algorithm
stops and returns \textbf{True} the condition
\begin{equation*}
\pi_u((f_c(G_1\times G_2)\cup LX)_I)\subset \mathbb R^u\setminus
\overline{B_u}(0,1)
\end{equation*}
is satisfied, which implies that for $t\in[0,1]$, $h(t,p,q)\notin
M_c.$\qed

Now we discuss how we verify condition (\ref{eq:hc3}). The main
point of our approach is that it is enough to compute $f_c(\bd
N_c)$.
\begin{alg}\label{alg:boundary}
\end{alg}
    \textbf{function} ComputeBoundary($\mathcal G$ : grid) : \textbf{bool}\\
    \textbf{var}\\
    \mbox{}\qquad $X$ : representable set;\\
    \textbf{begin}\\
    \mbox{}\qquad\textbf{foreach} $G\in\mathcal G$\\
    \mbox{}\qquad\textbf{begin}\\
        \mbox{}\qquad\qquad $X:=(f_c(G))_I$;\\
        \mbox{}\qquad\qquad\textbf{if not} $\pi_s(X)\subset\inte B_s(0,1)$ \textbf{return False};\\
    \mbox{}\qquad\textbf{end;}\\
    \mbox{}\qquad\textbf{return True;}\\
    \textbf{end.}
\begin{lemma}\label{lem:boundary}
Assume $N$, $M$ be a h-sets, $f:\mathbb R^n\to\mathbb R^n$ be such
that $f_c$ is a diffeomorphism. Let $\mathcal G$ be a grid of $\bd
N_c$. If Algorithm~\ref{alg:boundary} is called with argument
$\mathcal G$ and returns \textbf{\em True} then the homotopy
defined by (\ref{eq:convhom}) satisfies condition (\ref{eq:hc3}).
\end{lemma}
\textbf{Proof:} Let $x\in \bd N_c$. Since Algorithm stops and
returns \textbf{True} we obtain that $\pi_s(f_c(x))\in B_s(0,1)$.
From Eq.~(\ref{eq:convhom}) it follows that for $x\in\bd N_c$ and
$t\in[0,1]$ $\pi_s(h(t,x))\in B_s(0,1)$. Therefore for $x\in\bd
N_c$ and $t\in[0,1]$
\begin{equation}\label{eq:instable}
h(t,x)\notin \overline{B_u}(0,1)\times\bd B_s(0,1) = M_c^+.
\end{equation}
There remains to prove that Eq.~(\ref{eq:instable}) is satisfied
for $x\in N_c$ and $t\in[0,1]$. Since $f_c$ is a diffeomorphism,
the Brouwer--Jordan Theorem implies that $\pi_s(f_c(N_c))\subset
B_s(0,1)$. Hence, for $t\in[0,1]$
\begin{equation*}
\pi_s(h(t,N_c))\subset B_s(0,1)
\end{equation*}
which implies that $h(t,N_c)\cap M_c^+=\emptyset$.\qed

\subsection{Technical data.} The grids used in the numerical proof
of Lemma~\ref{lem:symcover} and Lemma~\ref{lem:covchain} always
consist of  ``boxes'', i.e.,  products of representable intervals.
The total number of boxes used in the proof is approximately
$2.2\cdot10^{8}$. The numerical proof of Lemma~\ref{lem:symcover}
and Lemma~\ref{lem:covchain} took approximately 36 minutes on
2.4GHz processor under PLD Linux Distribution.

The C++ sources with a short description how to run the program
are available at \cite{W1}.